\documentclass[letterpaper, 10 pt, conference]{ieeeconf}  

\IEEEoverridecommandlockouts                              

\usepackage{times} 
\usepackage{amsmath} 
\usepackage{amssymb}  

\usepackage{hyperref}
\hypersetup{colorlinks, urlcolor = blue, linkcolor = black}
\usepackage{tikz}
\usetikzlibrary{shapes,arrows,decorations.markings,positioning}
\usepackage{xspace}

\title{\LARGE \bf
An Accretive Operator Approach to Ergodic Problems for Zero-Sum Games
}

\author{Antoine Hochart$^{* \dag}$
\thanks{$^*$The author is with INRIA Saclay-Ile-de-France and CMAP, Ecole polytechnique, Route de Saclay, 91128 Palaiseau cedex, France.
{\tt\small antoine.hochart@cmap.polytechnique.fr}}
\thanks{$^\dag$The author is supported by a PhD fellowship of Fondation Math\'ematique Jacques Hadamard (FMJH)}%
}

\renewcommand{\theenumi}{\roman{enumi}}

\renewcommand{\geq}{\geqslant}
\renewcommand{\leq}{\leqslant}
\newcommand{\R}{\mathbb{R}}
\newcommand{\N}{\mathbb{N}}
\newcommand{\unit}{e}
\newcommand{\id}{Id}
\newcommand{\X}{\mathcal{X}}
\newcommand{\Y}{\mathcal{Y}}
\newcommand{\V}{\mathcal{V}}
\renewcommand{\S}{\mathcal{S}}
\newcommand{\state}{[n]}
\newcommand{\MIN}{{\small \textsf{MIN}}\xspace}
\newcommand{\MAX}{{\small \textsf{MAX}}\xspace}
\newcommand{\Hnorm}[1]{\| \ifx\\#1\\ \cdot \else #1 \fi \|_{\text{H}}}
\newcommand{\TP}{\mathbb{TP}}
\def\<#1,#2>{\langle #1, #2 \rangle}
\DeclareMathOperator{\dom}{dom}
\DeclareMathOperator{\range}{rg}
\DeclareMathOperator{\FP}{FP}

\newtheorem{theorem}{Theorem}
\newtheorem{lemma}[theorem]{Lemma}
\newtheorem{proposition}[theorem]{Proposition}
\newtheorem{corollary}[theorem]{Corollary}
\newtheorem{definition}[theorem]{Definition}
\newtheorem{example}[theorem]{Example}
\newtheorem{remark}[theorem]{Remark}

\begin{document}

\maketitle
\thispagestyle{empty}
\pagestyle{empty}

\begin{abstract}

Mean payoff stochastic games can be studied by means of a nonlinear spectral problem involving the Shapley operator: the ergodic equation.
A solution consists in a scalar, called the ergodic constant, and a vector, called bias.
The existence of such a pair entails that the mean payoff per time unit is equal to the ergodic constant for any initial state, and the bias gives stationary strategies.
By exploiting two fundamental properties of Shapley operators, monotonicity and additive homogeneity, we give a necessary and sufficient condition for the solvability of the ergodic equation for all the Shapley operators obtained by perturbation of the transition payments of a given stochastic game with finite state space.
If the latter condition is satisfied, we establish that the bias is unique (up to an additive constant) for a generic perturbation of the transition payments.
To show these results, we use the theory of accretive operators, and prove in particular some surjectivity condition.

\end{abstract}

\section{INTRODUCTION}

In this paper, we are interested in stochastic games with mean payoff.
A first question is to understand under which conditions the value of such games is independent of the initial state.
This question can be studied by means of the {\em ergodic equation}.
When the state space is finite, the latter writes
\begin{equation}
  \label{eq:ergodic-equation}
  T(u) = \lambda \unit + u \enspace ,
\end{equation}
where $T: \R^n \to \R^n$ is the Shapley operator, $\unit$ is the unit vector of $\R^n$, and $(\lambda,u) \in \R \times \R^n$ is a solution.
If such a pair exists, the scalar $\lambda$, called {\em ergodic constant}, gives the mean payoff per time unit for every initial state, and the vector $u$, called {\em bias}, determines optimal stationary strategies.
Here, we study the situation in which the ergodic equation has a solution for all the Shapley operators obtained by perturbation of the transition payments of a given stochastic game.
When this property is satisfied, we say that the game is {\em ergodic}.
Several results in control theory~\cite{KM97,AG03}, discrete event systems theory~\cite{Ols91,YZ04} or Perron-Frobenius theory~\cite{Nus88,GG04} leads to sufficient ergodicity conditions for games, which are related to accessibility conditions.
Under particular hypothesis (bounded transition payment function~\cite{AGH15a}, ``weakly convex'' Shapley operator~\cite{CCHH10}, etc.), these conditions are also necessary.

A second question concerns the structure of the set of bias.
In the one-player case (optimal control problems), characterizations have been given, both in the deterministic framework~\cite{BCOQ92} and the  stochastic one~\cite{AG03}.
The two-player case appears less accessible and understanding the situations in which the bias is unique (up to an additive constant) is already interesting.
In a previous work~\cite{AGH14}, we have shown the uniqueness of the bias for a generic perturbation of the payments in the case of ergodic games with perfect information and finite action spaces.

In the present work, we apply the theory of accretive operators to the study of zero-sum games.
We first establish a preliminary result on the surjectivity of accretive maps in finite dimension.
Then, we deduce some results concerning the fixed points of nonexpansive maps: we give a condition equivalent to the existence of a fixed point for all additive perturbations of such a given map, and we establish, under this condition, the uniqueness of the fixed point for a generic perturbation.
Finally, we apply these results to stochastic games with finite state space.
We show in particular that the sufficient ergodicity condition introduced in~\cite{AGH15b}, involving hypergraphs, is in fact necessary, and we solve a problem raised in~\cite{GG04}.

\section{PRELIMINARIES ON ACCRETIVE MAPS}

\subsection{Definitions}

Throughout the paper, $(\X, \|\cdot\|)$ is a finite-dimensional real vector space equipped with a given norm.
We denote by $A: \X \rightrightarrows \X$ a {\em set-valued map} $A$ from $\X$ to itself, i.e., a map from $\X$ to the powerset of $\X$.
The {\em domain} of $A$ is defined by $\dom (A) := \{ x \in \X \mid A(x) \neq \emptyset \}$, and its {\em range} by $\range (A) := \bigcup_{x \in \X} A(x)$.
The {\em inverse} of $A$, denoted by $A^{-1}$, is the set-valued map from $\X$ to $\X$ sending $y \in \X$ to $\{ x \in \X \mid y \in A(x) \}$, so that $x \in A^{-1}(y)$ if, and only if, $y \in A(x)$.
For background on set-valued maps we refer the reader to the monograph of Aubin and Frankowska~\cite{AF09}.

We denote by $\X^*$ the dual space of $\X$, by $\|\cdot\|^*$ its dual norm, and by $\<\cdot,\cdot>$ the duality product.
The {\em (normalized) duality mapping} on $\X$ is the set-valued map $J: \X \rightrightarrows \X^*$ defined by
\[
  J(x) := \{x^* \in \X^* \mid \|x^*\|^* = \|x\|, \enspace \<x,x^*> = \|x^*\|^* \: \|x\| \} \enspace .
\]
Note that $\dom(J) = \X$ by the Hahn-Banach theorem, and that $J(x)$ is a compact convex subset of $\X^*$ for every $x \in \X$.

A set-valued map $A: \X \rightrightarrows \X$ is {\em accretive} if, for every $x,y \in \X$, every $u \in A(x)$ and every $v \in A(y)$, there exists an element $x^* \in J(x-y)$ such that $\<u-v,x^*> \geq 0$.
Furthermore, denoting by $\id$ the identity map on $\X$, if $\range(\id + \lambda A) = \X$ for some (hence all) $\lambda > 0$, then $A$ is {\em $m$-accretive}.
Finally, $A$ is {\em coaccretive} if its inverse $A^{-1}$ is accretive. 

Accretive maps have been widely studied as infinitesimal generators of nonexpansive semigroups associated to evolution equations in Banach spaces, see~\cite{Kat67,CL71,Bro76, Rei76}.
They also naturally appear when considering trajectories defined by nonexpansive maps, see~\cite{Vig10} and the references therein.
Indeed, if a map $T: \X \to \X$ is nonexpansive (with respect to the norm of $\X$), meaning that $\|T(x)-T(y)\| \leq \|x-y\|$ for all $x, y \in \X$, then the operator $A := \id - T$ is $m$-accretive. 

\subsection{Surjectivity Conditions}

The main result of this section is a necessary surjectivity condition for accretive maps in finite-dimensional spaces.

\begin{theorem}
  \label{thm:surjectivity-condition}
  Let $A: \X \rightrightarrows \X$ be an accretive map on a finite-dimensional vector space $\X$.
  If $\range (A) = \X$, then the subset $\S_\alpha := \{x \in \X \mid \inf_{u \in A(x)} \|u\| \leq \alpha \}$ is bounded for every scalar $\alpha \geq 0$.
  Moreover, if $A$ is $m$-accretive, then the converse also holds true.
\end{theorem}

The second part of Theorem~\ref{thm:surjectivity-condition} is a special case of a corollary of Th.~3 in~\cite{KS80}.
As for the first part, it is deduced from the subsequent proposition.
The latter is a transposition to coaccretive maps, at least in finite dimension, of the known fact that an accretive map $A$ in a reflexive Banach space is {\em locally bounded} at any point $x$ in the interior of its domain, meaning that there exists a neighborhood $\V$ of $x$ such that $A(\V)$ is bounded, see \cite{FHK72}.

\begin{proposition}
  \label{prop:local-boundedness}
  Let $A: \X \rightrightarrows \X$ be a coaccretive map and let $x$ be a point in the interior of $\dom(A)$.
  Then $A$ is locally bounded at $x$.
\end{proposition}

\begin{example}[{\em discrete $p$-Laplacian}]
  Let $G=(V,E)$ be a finite connected undirected graph.
  Every edge $\{i,j\}$ is equipped with a weight $C_{i j} > 0$.
  For $p>1$, the {\em discrete $p$-Laplacian} is the map $L_p: \R^V \to \R^V$, whose coordinate $i \in V$ is defined by
  \begin{equation*}
    \left( L_p(v) \right)_i = \sum_{j : \{i,j\} \in E} C_{i j} (v_i - v_j) |C_{i j} (v_i - v_j)|^{p-2} \enspace .
  \end{equation*}

  Let $B \subset V$ be a nonempty subset different from $V$, and let $w \in \R^B$.
  We are interested in the following boundary value problem: given $g \in \R^{V \setminus B}$, find $v \in \R^V$ such that
  \begin{equation}
    \label{eq:Dirichlet-p-Laplacian}
    \begin{cases}
      \left( L_p(v) \right)_i = -g_i \enspace , \quad &\forall i \in V \setminus B \enspace ,\\
      v_i = w_i \enspace , \quad &\forall i \in B \enspace .
    \end{cases}
  \end{equation}
  Note that this problem is equivalent to the minimization of the energy function
  \[
    \sum_{\{i,j\} \in E} \frac{1}{p C_{i j}} |C_{i j} (v_i - v_j)|^p + \sum_{i \in V \setminus B} g_i v_i \enspace ,
  \]
  subject to $v_i = w_i$ for all $i \in B$.
  In particular, when $p=2$, it recovers the classical problem of computing the electrical potential $v$ on the graph $G$, with a prescribed potential $w_i$ at node $i \in B$ and a prescribed current $g_i$ at node $i \in V \setminus B$, $C_{i j}$ being the conductance of the edge $\{i,j\}$.

  Let us reformulate Problem~\eqref{eq:Dirichlet-p-Laplacian} as follows.
  Let $\X=\R^{V \setminus B}$ and, for $x \in \X$ and $w \in \R^B$, let $x|w$ be the vector in $\R^V$ given by $(x|w)_i = x_i$ for $i \in V \setminus B$ and $(x|w)_i = w_i$ for $i \in B$.
  Let us now introduce the operator $A: \X \to \X$ whose coordinate $i \in V \setminus B$ is defined by $\left( A(x) \right)_i = \left( L_p(x|w) \right)_i$.
  Thus, Problem~\eqref{eq:Dirichlet-p-Laplacian} is equivalent to finding a solution $x \in \X$ to the equation $A(x)=-g$.
  
  It can be shown that $A$ is $m$-accretive (it is the gradient of a convex continuous function).
  Hence, according to Theorem~\ref{thm:surjectivity-condition}, the latter equation has a solution for every $g \in \X$ if, and only if, the subsets $\S_\alpha = \{x \in \X \mid \|A(x)\| \leq \alpha \}$ are bounded for all $\alpha > 0$ ($\|\cdot\|$ being any norm).
  We leave it to the reader to check that this condition holds.
\end{example}

\section{FIXED POINT PROBLEMS FOR NONEXPANSIVE MAPS}

\subsection{Existence Stability of Fixed Points}

As a direct application of Theorem~\ref{thm:surjectivity-condition}, we get a necessary and sufficient condition for the stability, under additive perturbations, of the existence of fixed points of a nonexpansive operator.

\begin{corollary}
  \label{coro:stability-FP}
  Let $(\X, \|\cdot\|)$ be a finite-dimensional real vector space and let $T: \X \to \X$ be a nonexpansive operator.
  For every vector $g \in \X$, the operator $g+T$ has at least one fixed point if, and only if, for every scalar $\alpha \geq 0$, the set $\S_\alpha := \{x \in \X \mid \|x-T(x)\| \leq \alpha \}$ is bounded.
\end{corollary}

Note that the ``if'' part is also readily obtained as a corollary of Th.~4.1 in~\cite{Nus88}.

\subsection{Generic Uniqueness of the Fixed Point}

We now fix a nonexpansive operator $T: \X \to \X$ and assume that the condition in Corollary~\ref{coro:stability-FP} is satisfied.
Then, we define the set-valued map $\FP: \X \rightrightarrows \X$ by $\FP(g) := \{x \in \X \mid x - T(x) = g\}$, i.e., the map that associates to a vector $g \in \X$ the set of fixed points of the operator $g+T$.
Note that, by hypothesis, we have $\dom(\FP) = \X$ and $\FP^{-1} = \id - T$.

We first deduce from the condition in Corollary~\ref{coro:stability-FP} some properties on $\FP$.

\begin{lemma}
  \label{lem:FP-usc}
  The set-valued map $\FP$ has compact values and is upper semicontinuous.
\end{lemma} 

Then, by exploiting the accretivity of the operator $\id - T$ we get a characterization of the vectors $g$ for which $g+T$ has a unique fixed point.
\begin{theorem}
  \label{thm:uniqueness-continuity}
  The set-valued map $\FP$ is continuous at point $g \in \X$ if, and only if, $\FP(g)$ is a singleton.
\end{theorem}

The properties mentioned in Theorem~\ref{thm:uniqueness-continuity} are in fact {\em generic}, meaning that the set of elements of $\X$ for which $\FP$ is single-valued is a {\em residual}.
Recall that a {\em residual} of $\X$ is a countable intersection of dense open subsets of $\X$.
According to Baire's Theorem, a residual of the finite-dimensional vector space $\X$ is dense.
This generic property is stated in the following result.

\begin{theorem}[see {\cite[Th.~1.4.13]{AF09}}]
  \label{thm:generic-continuity-FP}
  Let $\X$ and $\Y$ be two finite-dimensional vector spaces, and let $A: \X \rightrightarrows \Y$ be an upper semicontinuous set-valued map.
  Then, $A$ is continuous on a residual of $\X$.
\end{theorem}

\section{APPLICATION TO GAMES WITH ERGODIC PAYOFF}

\subsection{Stochastic Games with Mean Payoff}

An important application where nonexpansive operators appears, which is our main motivation, is two-player zero-sum repeated games.
We consider here the case of stochastic games with a finite state space, say $\state := \{1,\dots,n\}$.
When the current state is $i \in \state$, the (mixed) action spaces of the Players \MIN and \MAX are denoted by $A_i$ and $B_i$, respectively, and if actions $a \in A_i$ and $b \in B_i$ are selected by the players, the transition payment is denoted by $r_i^{a b} \in \R$ and the transition probability by $P_i^{a b} := (P_{i j}^{a b})_{j \in \state} \in \R^n$.

The game is played in stages, starting from a given initial state $i_0$, as follows: at step $\ell$, if the current state is $i_\ell$, the players choose actions $a_\ell \in A_{i_\ell}$ and $b_\ell \in B_{i_\ell}$.
Then, Player \MIN pays $r_{i_\ell}^{a_\ell b_\ell}$ to Player \MAX and the next state is chosen according to the probability law $P_{i_\ell}^{a_\ell b_\ell}$.

We consider here games played in finite horizon.
The payoff of a play in the $k$-stage game (played in horizon $k$) is given by the sum of the transition payments of the $k$ first stages, namely $\Sigma_{0 \leq \ell \leq k-1} r_{i_\ell}^{a_\ell b_\ell}$.
Player \MAX intends to maximize this quantity while Player \MIN intends to minimize it.
This gives rise to the value (if it exists) of the $k$-stage game, which we denote by $v_i^k$ when the initial state is $i$.

The value vector $v^k = (v_i^k)_{i \in \state} \in \R^n$, which we assume to exist, satisfies the following recursive formula:
\begin{equation}
  \label{eq:recursive-structure}
  v^k = T(v^{k-1}) \enspace , \quad v^0 = 0 \enspace ,
\end{equation}
where $T: \R^n \to \R^n$ is the Shapley operator, whose $i$th coordinate is given by
\begin{equation}
  \label{eq:Shapley-operator}
  T_i(x) = \inf_{a \in A_i} \sup_{b \in B_i} \Big( r_i^{a b} + \sum_{j \in \state} P_{i j}^{a b} x_j \Big) \enspace .
\end{equation}
Here, $\inf$ and $\sup$ commutes.
The asymptotic behavior of the average payoff vector per time unit, $v^k / k$, as the number of stages $k$ grows to infinity is a major topic in game theory, see~\cite{NS03}.
When the limit exist, we call it the {\em mean-payoff vector}.
Because of the recursive structure of the game expressed in~\eqref{eq:recursive-structure}, this problem amounts to study the orbit $\{ T^k(0) \mid k \in \N \}$ of the Shapley operator.
In particular, note that if the nonlinear spectral problem~\eqref{eq:ergodic-equation} has a solution, then the sequence $(v_i^k / k)$ converges to $\lambda$ for any state $i \in \state$.
This fact relies on the following two properties satisfied by any Shapley operator $T$:
\begin{itemize}
  \item monotonicity: $x \leq y \implies T(x) \leq T(y)$, \, $x, y \in \R^n$;
  \item additive homogeneity: $T(x + \lambda \unit) = T(x) + \lambda \unit$, \, $x \in \R^n$, \, $\lambda \in \R$.
\end{itemize}
This implies in particular that $T$ is sup-norm nonexpansive.

Let us mention that any operator on $\R^n$ that is monotone and additively homogeneous can be written as a Shapley operator~\eqref{eq:Shapley-operator}, see~\cite{Kol92,RS01}.
This motivates the following definition.
\begin{definition}
  We call {\em Shapley operator} (on $\R^n$) any operator from $\R^n$ to itself that is both monotone and additively homogeneous.
\end{definition}

\subsection{Ergodicity Conditions for Zero-Sum Games}

Let us introduce the {\em Hilbert's seminorm} on $\R^n$, defined by
\[
  \Hnorm{x} = \max_{1 \leq i \leq n} x_i - \min_{1 \leq i \leq n} x_i \enspace .
\]
Let $\TP^n$ be the ``additive projective space'' of $\R^n$, defined as the set of equivalence classes of $\R^n$ where two vectors $x, y \in \R^n$ are equivalent if there exists a scalar $\alpha \in \R$ such that $x-y = \alpha \unit$.
Note that, in this case, we have $\Hnorm{x} = \Hnorm{y}$.
Hence, the space $(\TP^n, \Hnorm{})$ is a normed space with finite dimension.

It is a standard result (see for instance~\cite{GG04}) that any Shapley operator $T: \R^n \to \R^n$ is nonexpansive with respect to the Hilbert's seminorm.
Moreover, since $T$ is monotone and additively homogeneous, it can be quotiented into a map $[T]: \TP^n \to \TP^n$.
Thus, this quotiented map is nonexpansive with respect to the norm $\Hnorm{}$.

Observe now that a vector $u \in \R^n$ is a bias of $T$, meaning that there exists a scalar $\lambda \in \R$ such that the pair $(\lambda, u)$ solves the ergodic equation~\eqref{eq:ergodic-equation}, if, and only if, its equivalence class $[u] \in \TP^n$ is a fixed point of $[T]$.

Let $T$ be a Shapley operator on $\R^n$.
Given two scalars $\alpha, \beta \in \R$, we define the {\em slice space} $\S_\alpha^\beta$ as the subset of $\R^n$ $\{ x \in \R^n \mid \alpha \unit + x \leq T(x) \leq \beta \unit + x \}$. 
As a special case of Th.~4.1 in~\cite{Nus88} applied to the quotiented map $[T]$, we deduce that if all the slice spaces are bounded in the Hilbert's seminorm, then the game associated to the Shapley operator $T$ is {\em ergodic}, meaning that for every vector $g \in \R^n$, there is a scalar $\lambda \in \R$ and a vector $u \in \R^n$ such that $g+T(u) = \lambda \unit + u$.
Observe that, if $T$ is written as in~\eqref{eq:Shapley-operator}, then $g+T$ is the Shapley operator of a game almost identical to the one associated with $T$, except that the transition payments are given by $g_i + r_i^{a b}$, i.e., are perturbed by quantities that only depend on the state.

The above remarks and a direct application of Corollary~\ref{coro:stability-FP} to the quotiented map $[T]$ leads to necessary and sufficient ergodicity conditions for zero-sum repeated games.

\begin{theorem}
  \label{thm:ergodicity-conditions}
  Let $T: \R^n \to \R^n$ be a Shapley operator.
  The following are equivalent:
  \begin{enumerate}
    \item for every vector $g \in \R^n$, there is a pair $(\lambda,u) \in \R \times \R^n$ such that $g+T(u) = \lambda \unit + u$;
    \item the subset $\S_\alpha = \{ x \in \R^n \mid \Hnorm{x-T(x)} \leq \alpha \}$ is bounded in the Hilbert's seminorm for all $\alpha \in \R$;
        \item the slice space $\S_\alpha^\beta = \{ x \in \R^n \mid \alpha \unit + x \leq T(x) \leq \beta \unit + x \}$ is bounded in the Hilbert's seminorm for all $\alpha, \beta \in \R$.
  \end{enumerate}
\end{theorem}

\begin{example}
  \label{ex:ergodicity}
  Consider the Shapley operator $T: \R^3 \to \R^3$ given by
  \small
  \begin{equation}
    \label{eq:Shapley-example}
    T(x) =
    \begin{pmatrix}
      \sup_{0 < p \leq 1} \left( \; \log p + p ( x_1 \wedge x_3 ) + (1-p) x_2 \; \right) \\
      \inf_{0 < p \leq 1} \left( \; -\log p + p ( x_2 \vee x_3 ) + (1-p) x_1 \; \right) \\
      x_3
    \end{pmatrix}
  \end{equation}
  \normalsize
  where $\wedge$ and $\vee$ stand for $\min$ and $\max$, respectively.
  It corresponds to a game with three states.
  Player \MAX controls state~$1$, Player \MIN controls state~$2$ and state~$3$ is an absorbing state.
  In state~$1$, Player \MAX chooses an action $p \in (0,1]$ and receives $\log p$ from Player \MIN.
  Then, with probability $p$, Player \MIN selects the next state between $1$ and $3$, and with probability $1-p$ the next state is $2$.
  Hence, in state~$1$, Player \MAX is losing and she can only increase the probability to leave the state to a more profitable one (namely, state~$2$) at the expense of a blowing cost.
  A dual trade-off applies to Player \MIN in state~$2$.

  Using the hypergraph conditions in~\cite{AGH15b} for instance, we can show that all the slice spaces are bounded in the Hilbert's seminorm.
  Thus, according to Theorem~\ref{thm:ergodicity-conditions}, the value of the game does not depend on the initial state for any additional payment $g_i$ granted in state $i$.
\end{example}

\subsection{Generic Uniqueness of the Bias of Zero-Sum Games}

We now fix a Shapley operator $T: \R^n \to R^n$ and we assume that the associated game is ergodic.
A direct application of Lemma~\ref{lem:FP-usc}, Theorem~\ref{thm:uniqueness-continuity} and Theorem~\ref{thm:generic-continuity-FP} to the quotiented map $[T]$ leads to the following.

\begin{theorem}
  \label{thm:generic-uniqueness-bias}
  The set of perturbation vectors $g \in \R^n$ for which $g+T$ has a unique bias (up to an additive constant) is a residual of $\R^n$.
\end{theorem}

\begin{example}
  Consider the Shapley operator~\eqref{eq:Shapley-example} introduced in Example~\ref{ex:ergodicity}.
  We have already noted that the game associated to this operator is ergodic, and it is readily seen that, for a perturbation vector $g \in \R^3$, the ergodic constant of $g+T$ is $g_3$.
  Thus, for a generic vector $g \in \R^3$, there is a unique vector $u \in \R^3$ (up to an additive constant) such that $g+T(u) = g_3 \unit + u$.
\end{example}

\begin{remark}
  Let us mention that when $T$ is polyhedral, that is, when considering a game with {\em perfect information} and finite action spaces, we have more information on the geometry of the residual set.
  More specifically, we have shown in~\cite{AGH14} that the set of perturbation vectors $g$ for which $g+T$ has more than one bias is included in the finite union of subspaces of codimension $1$.
  
  As an example, let us consider the following Shapley operator defined on $\R^3$ and already introduced in~\cite{AGH14} (recall that the addition has precedence over $\min$ and $\max$, denoted here by $\wedge$ and $\vee$, respectively):
  \small
  \begin{equation*}
    T(x) =
    \begin{pmatrix}
      \frac{1}{2} (x_1 + x_3) \, \wedge \, 1 + \frac{1}{2} (x_1 + x_2)\\
      2 + \frac{1}{2} (x_1 + x_3) \, \wedge \, \left(1 + \frac{1}{2} (x_1 + x_2) \, \vee \, -2 + x_3 \right)\\
      3 + \frac{1}{2} (x_1 + x_3) \, \vee \, 1 + x_3
    \end{pmatrix} \enspace .
  \end{equation*}
  \normalsize
  It can be proved that the ergodic equation~\eqref{eq:ergodic-equation} is solvable for every perturbation vector $g \in \R^3$ (see~\cite{AGH15a} for an effective method).
  Figure~\ref{fig} shows the intersection between the hyperplane $\{g \in \R^3 \mid g_3 = 0\}$ and the locus of perturbation vectors $g$ for which $g+T$ may have more than one bias.
  Hence, for each $g$ in the interior of a full-dimensional polyhedron, $g+T$ has a unique bias (up to an additive constant).
  \begin{figure}[!h]
    \centering
    \begin{tikzpicture}[scale=0.13]
      \draw [->, >= angle 60,black] (-5,0) -- (20,0);0
      \draw (19,0) node[above] {$g_1$};
      \foreach \x in {0,5,10,15}
      \draw [black] (\x,0) -- (\x,-0.5);
      \draw (10,0) node[below] {\footnotesize $10$};
      \draw [->, >= angle 60,black] (0,-20) -- (0,7);
      \draw (0,6) node[right] {$g_2$};
      \foreach \y in {-15,-10,-5,0,5}
      \draw [black] (-0.5,\y) -- (0,\y);
      \draw (0,-10) node[left] {\footnotesize $-10$};
      \draw (15,6) node {$g_3=0$};
      \draw [thick,blue] (5,-20) -- (5,7);
      \draw [thick,blue] (-2,7) -- (23/2,-20);
      \draw [thick,blue] (11,-11) -- (-5,5);
      \draw [thick,blue] (11,-11) -- (20,-42/3);
      \draw [thick,blue] (11,-11) -- (16,-20);
      \draw (0,0) node {$\bullet$};
      \draw (0,0) node [below left] {\footnotesize $0$};
    \end{tikzpicture}
    \caption{}
    \label{fig}
  \end{figure}
\end{remark}



\bibliographystyle{IEEEtran}
\bibliography{IEEEabrv,references}

\end{document}